\documentstyle[fleqn,12pt]{article}
\begin{document}\pagenumbering{arabic}\setcounter{page}{1}\pagestyle{plain}
\baselineskip=16pt

\begin{center}{\large\bf On the determination of the relation between the 
matrix elements of CP$_q(n)$ and covariant $q$-oscillators} \end{center}

\vspace{2cm}
Salih CELIK \\
Department of Mathematics, Faculty of Sciences, Mimar Sinan University 
80690 Besiktas, Istanbul, TURKEY.

\vspace{3cm}
{\bf Abstract}:
In this article it is explicitly shown that the matrix elements 
of $CP_q(n)$ are the annihilation operators of Pusz-Woronowicz (P-W)
oscillators provided they are rescaled.

\vfill\eject
{\bf 1. Introduction }

Recently, using the first order differential calculus and the 
second quantization procedure the algebra of $q$-oscillators covariant 
under the coaction of the unitary quantum group $SU_q(n)$ was given by 
Pusz and Woronowicz [1]. Wess and Zumino [2] have developed a differential 
calculus for the quantum hyperplane which is covariant under 
the action of $GL_q(n)$. In Ref. 3 it has been shown that the unitary 
quantum group $U_q(n)$ can be constructed in terms of n(n-1)/2 
$q$-oscillators and n-commuting phases. This approach shows that the 
quantum projective space $CP_q(n)$ can be identified with n 
$q$-oscillators.

A new realization of the quantum group $SU_q(2)$ has been obtained 
by introducing a $q$-analogue of the usual harmonic oscillator and 
using the Jordan-Schwinger mapping [4,5]. For early studies on the 
$q$-oscillator see Arik et al [6,7]. Similar deformations of the 
quantum harmonic oscillator algebra have attracted a lot of 
attention [8-17].

The purpose of this study is to show that the elements of the 
matrix representations of the quantum projective space $CP_q(n)$, 
by using the representations found in [3], are the covariant 
annihilationoperators of Pusz-Woronowicz (P-W) oscillators.

In section 2 we introduced the quantum group $SU_q(2)$. And also 
we showed the matrix elements of $SU_q(2)$ can map to a Biedenharn-Macfarlane 
(B-M) oscillator algebra [4,5]. In section 3 we proved that the matrix 
elements of $CP_q(n)$ are the annihilation operators of (P-W) oscillators. 
The results has been discussed in section 4.

{\bf 2. The quantum group $SU_q(2)$ and $q$-oscillators }

The quantum group $SU_q(2)$ consists of all matrices of the form
$$ M = \left(\matrix{ a  &  b \cr
                    -qb^* & a^* \cr }
\right) $$
where
$$ ab = qba, ~~ab^* = qb^*a, ~~bb^* = b^*b, \eqno(2.1) $$
$$ aa^* + q^2bb^* = 1, ~~ a^*a + bb^* = 1, \eqno(2.2) $$
and $0 < q < 1 ~(q \in {\cal R})$.

The $^*$-operation is defined as follows:
$$ (a)^* = a^+  ~~\mbox{and}~~  (a^+)^* = a. $$
From the equations in (2.2) one obtains
$$ aa^* - q^2a^*a = 1 - q^2. \eqno(2.3) $$
This equation can be used to define our $q$-oscillator [3] where the 
operator $a^*$ is the creation and its hermitean conjugate $a$ is 
the annihilation operator.Note that the normalization of the 
oscillators is chosen so that in the $q \rightarrow 1$ limit this 
operators commute and gives c-numbers. The spectrum of the 
$q$-oscillator defined by equation (2.3) is completely fixed for 
$0 < q < 1$ and it is shown that [3]
$$ a^*a = 1 - q^{2N} = [N]_q \eqno(2.4) $$
where $N$ is the number operator. From this we can get immediately
$$ [N,a] = - a, ~~[N,a^*] = a^*. \eqno(2.5) $$

Now if $a$ is defined as follows
$$ a = (1 - q^2)^{1/2}~q^{N/2} ~\mbox{a}_q $$
then
$$ a^* = (1 - q^2)^{1/2}~\mbox{a}_q^+ ~q^{N/2} $$
is obtained. Here $N^+ = N$. Then (2.3) becomes
$$ \mbox{a}_q \mbox{a}_q^+ - q \mbox{a}_q^+ \mbox{a}_q = q^{-N}. \eqno(2.6) $$
It follows that the matrix elements of the two-dimensional representation of 
$SU_q(2)$ can be mapped to a (B-M) oscillator algebra,
i.e. the algebra (2.3,5) is equivalent to (B-M) oscillator algebra. 
Note that although for $0 < q < 1$ equations (2.3) and (2.6) are related 
by just a rescaling, in the $q \rightarrow 1$ limit they are 
fundamentally different. In this limit the oscillator creation and 
annihilation operators in (2.3) are commutative so that $a$ and 
$a^*$ can be considered as commuting complex numbers.

{\bf 3. The quantum projective space $CP_q(n)$ 
and covariant $q$-oscillators}

Now the quantum group $SU_q(n)$ is discussed. It can be shown 
that any element of $SU_q(n)$ can be expressed uniquely as [3]
$$ M = \prod_{k=1}^{n-1} M_{i,i+1}(a_{ik}) \prod_{i=1}^{n-1} \chi_{i,i+1}
(\beta_i) $$
and 
$$ H = \prod_{i=1}^{n}M_{i,i+1}(a_i) \eqno(3.1) $$
being identified with the matrix representatives of the n-dimensional 
quantum projective space $CP_q(n)$. Here each $M_{i,i+1}(a_i)$ is 
the matrix whose 2x2 diagonal block in the $i,i+1$ position is a 
$q$-oscillator matrix
$$ \left( \matrix{  a_i & (1 - a_i^* a_i)^{1/2} \cr
               -q(1 - a_i^* a_i)^{1/2} & a_i^* \cr}
\right) $$
and the remaining diagonal elements contain 1, and all the elements 
apart from these are zero. In (3.1) each $a_i$ is a $q$-oscillator 
and their properties are determined by the commutation relations,
$$ [a_i,a_j] = 0 = [a_i^*,a_j^*], ~~i,j=1,2,...,n \eqno(3.2a) $$
$$ [a_i,a_j^*] = (1 - q^2) q^{2N_i} \delta_{ij}, \eqno(3.2b) $$
$$ [a_i,N_j] = \delta_{ij}a_i, ~~[a_i^*,N_j] = - \delta_{ij}a_i^*, 
\eqno(3.2c) $$
where $\delta_{ij}$ denotes the Kronecker delta.

In the matrix $H$ given by (3.1) instead of taking the n-1 
$a_i$'s as independent one can take the first n-1 elements of the 
first row as independent. Then it can be shown that these 
elements can be shown  they are the annihilation operators of 
(P-W) oscillators [1]. Note that the first elements taken are
$$ a = (1 - q^2)^{1/2}~A $$
$$ a^* = (1 - q^2)^{1/2}~A^+. \eqno(3.3) $$

In fact, if we take
$$ A_k = (1 - q^2)^{-1/2} \left(\prod_{i=1}^{k-1}x_i \right)a_k, 
~~\mbox{for} k = 2,3,...,n \eqno(3.4) $$
where $x_i = (1 - a_i^*a_i)^{1/2}$ and 
$$ A_1 = (1 - q^2)^{-1/2}a_1 $$
one has
$$ A_iA_j = qA_jA_i, ~~i < j $$
$$ A_iA_j^+ = qA_j^+A_i, ~~i \neq j  \eqno(3.5) $$
$$ A_i^+A_j^+ = qA_j^+A_i^+, i > j $$
$$ A_iA_i^+ - q^2A_i^+A_i = 1 + (q^2 - 1) \sum_{k<i}A_k^+4A_k. $$

Now let's take
$$ a_1 = A_1. $$
For a representation in a Hilbert space where $a_1^*$ denotes the 
hermitean conjugate of $a_1$, the eigenvalues of $A_1^+A_1$ are real 
and non-negative. Moreover, by the relation (2.2) the operator 
$A_1^+A_1$ is bounded and the eigenvalues of this operator are in 
the interval [0,1). Thus one can take
$$ a_2 = (1 - A_1^+A_1)^{-1/2}A_2. $$
We require that the relations which will be obtained, after the 
defining the transformations, are the same with the relations (3.2).
Thereby we will take
$$ a_1 = (1 - q^2)^{1/2}A_1, \eqno(3.6a) $$
$$ a_2 = (1 - q^2)^{1/2}X^{-1}A_2, \eqno(3.6b) $$
where
$$ X^2 = A_1A_1^+ - A_1^+A_1. $$
Note that
$$ A_1X^2 = q^2X^2A_1,  ~~A_1^+X^2 = q^{-2}X^2A_1^+, $$
$$ A_1X^{\pm 1} = q^{\pm 1}X^{\pm 1}A_1, ~~
   X^{\pm 1}A_1^+ = q^{\pm 1}X^{\pm 1}A_1^+, \eqno(3.7) $$
$$ [A_2,X^2] = 0 = [A_2^+,X^2]. $$
Now by using the commutation relations (3.5) and (3.7) it is easy to 
check that the relations (3.2) are invariant under the transformations (3.6). 
Indeed, for example,
$$ a_1a_2= (1 - q^2)A_1X^{-1}A_2 = (1 - q^2)X^{-1}A_2A_1 = a_2a_1. $$
Consequently, we have shown that there is a one to one correspondence 
between the matrix elements of the matrix representation of 
$CP_q(n)$ and covariant (P-W) oscillators.

In the general case, one can take 
$$ a_1 = (1 - q^2)^{1/2}A_1 $$
and for $k = 2,3,...,n $
$$ a_k=(1 - q^2)^{1/2} \left(\prod_{i=1}^{k-1}X_i^{-1} \right)A_k \eqno(3.8) $$
where 
$$ X_i^2 = (A_iA_i^+ - A_i^+A_i) $$
and a straightforward computation shows that the relations (3.2) are 
satisfied.

{\bf 4. Discussion}

In this study by taking the matrix elements of the unitary 
quantum group as the independent $q$-oscillators we proved that the 
entries of the matrix representations of the quantum projective 
space $CP_q(n)$ are the annihilation operators of Pusz-Woronowicz 
oscillators.

The quantum projective space has two different 
$q \rightarrow 1$ limits. In the first case it reduces to the 
ordinary complex projective space $CP(n)$ and in the other case it 
reduces to the n-dimensional ordinary oscillator.

{\bf References}

\noindent
$[1]$ W. Pusz and S. L. Woronowicz, Rep. Math. Phys. 27 (1989) 231.\\
$[2]$ J. Wess and B. Zumino, Nucl. Phys. B (Proc.Suppl) 18 (1990) 302.\\
$[3]$ M. Arik and S. Celik, Z.Phys. C. 59 (1993) 99.\\
$[4]$ A. J. Macfarlane, J. Phys. A : Math. Gen. 22 (1989) 4581.\\
$[5]$ L. C. Biedenharn, J. Phys. A : Math. Gen. 22 (1989) L873.\\
$[6]$ M. Arik, D. D. Coon and Y.Lam, J. Math. Phys. 16 (1975) 1765.\\
$[7]$ M. Arik and D. D. Coon, J. Math. Phys. 17 (1976) 524.\\
$[8]$ Y. J. Ng, J. Phys. A : Math. Gen. 23 (1990) 1023.\\
$[9]$ M. Chaichian, P. Kulish and J. Lukierski, Phys.Lett. B 237 (1990)401.\\
$[10]$ B. Fairlie and C. K. Zachos, Phys.Lett. B 256(1) (1991) 43.\\
$[11]$ M. Arik, Z.Phys. C 51 (1991) 627.\\
$[12]$ S. P. Vokos, J. Math. Phys. 32 (11) (1991) 2979. \\
$[13]$ M. Chaichicn, P.Kulish and J.Lukierski, Phys.Lett. B 262 (1991)43.\\
$[14]$ C. Daskaloyannis, J. Phys. A : Math. Gen. 23 (1990) L789.\\
$[15]$ P. Kulish and E. V. Damaskinsky, J.Phys.A:Math.Gen. 23 (1990) L415.\\
$[16]$ R. M. Mir-Kasimov, J.Phys. A: Math. Gen. 24 (1991) 4283.\\
$[17]$ M. Arik and M. Mungan, Phys. Lett. B 282 (1992) 101.
\end{document}